\documentclass[onefignum,onetabnum]{siamonline250211}
\usepackage{graphicx}
\usepackage{xcolor}
\usepackage{amssymb,amsfonts}
\usepackage{circuitikz}
\usetikzlibrary{patterns}
\usepackage{booktabs}

\newcommand{\R}{\mathbb R}
\newcommand{\N}{\mathbb N}
\newcommand{\C}{\mathbb C}
\newcommand{\intv}{\mathcal I}
\newcommand{\diff}{\mathrm d}
\newcommand{\lag}{\mathfrak l}
\newcommand{\calE}{\mathcal E}
\newcommand{\calF}{\mathcal F}
\newcommand{\calO}{\mathcal O}
\newcommand{\calT}{\mathcal T}
\newcommand{\calS}{\mathcal S}
\newcommand{\calK}{\mathcal K}
\newcommand{\calL}{\mathcal L}

\newcommand{\calP}{\mathcal P}

\newcommand{\calR}{\mathcal R}

\newcommand{\calD}{\mathcal D}
\newcommand{\calM}{\mathcal M}
\newcommand{\B}{\mathfrak B}
\DeclareMathOperator{\col}{col}

\let\ker\relax
\DeclareMathOperator{\ker}{ker}
\DeclareMathOperator{\rank}{rank}

\newsiamremark{remark}{Remark}
\crefname{remark}{remark}{remarks}
\Crefname{remark}{Remark}{Remarks}

\title{From data pencils to image representations of continuous-time LTI systems}
\author{Timo Reis\thanks{Systems Theory and Partial Differential Equations Group, Institute of Mathematics, Technische Universit\"at Ilmenau, 98693 Ilmenau, Germany (\email{timo.reis@tu-ilmenau.de}).}
\and Philipp Schmitz\thanks{Optimization-based Control Group, Institute of Mathematics, Technische Universit\"at Ilmenau, 98693 Ilmenau, Germany (\email{philipp.schmitz@tu-ilmenau.de}).}
\and Karl Worthmann\thanks{Optimization-based Control Group, Institute of Mathematics, Technische Universit\"at Ilmenau, 98693 Ilmenau, Germany (\email{karl.worthmann@tu-ilmenau.de}).}}
\headers{From data pencils to image representations}{T. Reis, P. Schmitz, and K. Worthmann}

\begin{document}

\maketitle

\begin{abstract}
    The continuous-time fundamental lemma describes the set of all trajectories of a linear time-invariant system by a Gramian matrix built from measured input--output data and their derivatives. The jet of a trajectory is represented as the Gramian applied
to a nonunique parameter satisfying a differential-algebraic equation. We determine the complete Kronecker structure of the associated
data pencil for multiple-input multiple-output systems. This enables the construction of a polynomial image representation of the
input--output behavior directly from the data. The resulting
image matrix has column degrees equal to the controllability
indices, and its latent variable is free and serves as a
flat output of the system. The computation relies on an orthogonal staircase reduction followed by basis reconstruction using polynomial kernels. We also propose a structure-preserving preconditioning of the Gramian to improve its condition number and facilitate rank decisions in the reduction. We illustrate our approach on a chain of mass-spring-damper systems.
\end{abstract}

\begin{keywords}
continuous-time fundamental lemma, data-driven systems, image representation,
Kronecker structure, minimal polynomial bases, staircase algorithm
\end{keywords}

\begin{MSCcodes}
93B05, 93B20, 15A22, 65F45
\end{MSCcodes}

\section{Introduction}
Data-driven control aims to design controllers directly from measured
trajectories, without first identifying a model. For controllable
discrete-time linear time-invariant (LTI) systems, the fundamental
lemma of Willems et al.~\cite{willemsNotePersistencyExcitation2005}
states that, under persistency of excitation, finite-length trajectories
are spanned by the columns of a data Hankel matrix. This result
underlies data-driven prediction and control
\cite{coulsonDataEnabledPredictiveControl2019,berberichDataDrivenModelPredictive2021,depersisFormulasDataDrivenControl2020};
see also the review papers~\cite{markovskyBehavioralSystemsTheory2021,faulwasserBehavioralTheoryStochastic2023}.

Related developments include data informativity
\cite{vanwaardeDataInformativity2020}, with applications to switched systems \cite{schmitzDataInformativityFeedback2026}, extensions to LPV and
stochastic systems
\cite{verhoekFundamentalLemmaLPV2021,faulwasserBehavioralTheoryStochastic2023},
and nonlinear generalizations, including kernel-based approaches
\cite{berberichLinearTrackingMPC2022,alsaltiDataBasedSystemAnalysis2021,lianNonlinearDataEnabledPrediction2021,molodchykExploringLinksFundamental2024}.

Many physical systems are naturally modeled in continuous time. Sampling turns them into discrete-time systems, but the resulting model depends on the sampling period and ignores the intersample behavior. This has motivated several continuous-time versions of the fundamental lemma. In~\cite{lopezContinuousTimeVersionWillems2022,lopezInputOutputContinuousTimeVersion2024}, the data and their time shifts form a time-varying Hankel matrix, and trajectories are generated by a time-varying parameter vector that solves a differential equation. The approach in~\cite{lopezContinuousTimeVersionWillems2022} requires state measurements, whereas~\cite{lopezInputOutputContinuousTimeVersion2024} uses input--output data together with their derivatives. In~\cite{rapisardaPersistencyExcitationCondition2023}, persistency of excitation is defined via the input and its derivatives, and it is shown to guarantee identifiability. The fundamental lemma in~\cite{rapisardaFundamentalLemmaContinuoustime2023} is based on Chebyshev polynomial expansions of a trajectory and its derivatives. Linear combinations of the resulting coefficient sequences generate the Chebyshev coefficients of every trajectory. In practice, the expansions are truncated, so the resulting trajectories are approximate. The approach in~\cite{schmitzContinuoustimeFundamentalLemma2024} builds a Gramian matrix from the input--output data and their derivatives. This yields an exact characterization, in which the parameter generating a trajectory satisfies a differential-algebraic equation (DAE). All these approaches need derivatives of the measured signals, either explicitly or through a polynomial expansion. The approach in~\cite{othmaneDataBasedImageRepresentation2026}
uses algebraic differentiators to estimate the required
derivatives from measured data. In contrast,
\cite{wolskiAnotherContinuoustimeFundamental2026}
employs modulating functions, transferring differentiation
to known test functions and thereby avoiding differentiation
of the measured signals. Related continuous-time developments concern informativity based
on linear functionals \cite{ohtaDataInformativityContinuousTime2024},
derivative-free input--state formulations
\cite{wakaikiDataDrivenControlContinuousTime2026},
and extensions of the Gramian approach to LPV systems
\cite{schmitzContinuousTimeGeneralizationLPV2026}.
Applications include optimal, predictive, and iterative learning
control
\cite{schmitzContinuoustimeFundamentalLemma2024,schmitzDatadrivenContinuoustimeOptimal2026,wolskiDataDrivenModelPredictive2024,chuDataDrivenIterativeLearning2026}.

In this paper, we build on the Gramian-based fundamental lemma of \cite{schmitzContinuoustimeFundamentalLemma2024}. There, the jet of an input--output trajectory is represented
as $\Gamma g$ with the data Gramian $\Gamma$. The parameter $g$ is not unique and is constrained by a DAE with a rectangular, hence singular, matrix pencil. This makes the parameterization hard to use in computations. The approach in~\cite{othmaneDataBasedImageRepresentation2026}
obtains an image representation by embedding a staircase-reduced
data pencil into a unimodular matrix, but its analysis is
restricted to single-output systems. For multiple outputs,
the data pencil may have a nontrivial common left kernel,
which obstructs a direct unimodular construction.
First, we determine the complete Kronecker structure of
the data pencil and show that its nontrivial part is
strictly equivalent to an augmented controllability pencil. Second, we exploit this structure in an orthogonal
staircase reduction to reconstruct a minimal polynomial
basis of the kernel. This yields an image representation
of the full input--output behavior
(\cref{cor:full-image}). The resulting image
matrix has column degrees equal to the controllability
indices of the underlying LTI system, and its latent
variable is a flat output. Third, we propose a scaling of the Gramian designed to improve its conditioning and facilitate rank decisions in the staircase reduction.

The paper is organized as follows. \Cref{sec:FL} recalls persistency of excitation and the fundamental lemma and states the problem. \Cref{sec:pencil} analyzes the Kronecker structure of the data pencil. \Cref{sec:staircase} presents the staircase reduction and the data-driven image representation. \Cref{sec:numerics} discusses preconditioning and illustrates the method on a chain of mass-spring-damper systems.

\paragraph{Notation} For matrices $A_1,\dots,A_n$ with the same number of columns, we write \[\col(A_1,\dots,A_n) := \begin{bmatrix}
A_1^\top & \cdots & A_n^\top
\end{bmatrix}{}^\top;\] this notation extends to matrix-valued functions. For an open interval $\intv$ and $k\in\N$ let $H^k(\intv, \R^d)$ be the $k$th order Sobolev space based on $L^2(\intv, \R^d)$, with $H^0(\intv,\R^d)=L^2(\intv,\R^d)$. For a positive integer $k$ and $f\in H^{k-1}(\intv,\R^d)$ we denote the jet of $f$ of length $k$ by \[\Lambda_{k}(f) :=\col(f, f^{(1)}, \dots, f^{(k-1)})\colon\intv\to\R^{kd}.\] For two functions $f\in H^{k-1}(\intv,\R^{d})$, $g\in H^{j-1}(\intv, \R^{e})$ we stack the individual jets, $\Lambda_{k,j}(f, g):= \col(\Lambda_{k}(f), \Lambda_{j}(g))\colon\intv\to\R^{kd+je}$.
We write $\R[s]$ and $\C[s]$ for polynomial rings,
and $\R(s)$ and $\C(s)$ for the corresponding fields
of rational functions. Polynomial differential operators
acting on nonsmooth functions are understood in the
distributional sense.

\section{Fundamental lemma and problem statement}
\label{sec:FL}

We recall persistency of excitation and the continuous-time
fundamental lemma from
\cite{schmitzContinuoustimeFundamentalLemma2024,schmitzDatadrivenContinuoustimeOptimal2026}.
The resulting parameterization leads to a differential-algebraic
equation that will be analyzed below.

We consider the LTI system
\begin{equation}
\label{eq:sys}
    \dot{x}(t) = Ax(t) + Bu(t) ,\quad y(t) = Cx(t) + Du(t)
\end{equation}
with $A\in\R^{n\times n}$, $B\in\R^{n\times m}$, $C\in\R^{p\times n}$, $D\in\R^{p\times m}$ on an open and finite time interval $\intv$.
The object of interest is the set of input--output trajectories of system~\eqref{eq:sys}, called the (\emph{manifest}) \emph{behavior}, given by
\[
    \B = \left\{(y,u)\in L^2(\intv, \R^p)\times L^2(\intv, \R^m) \,\middle|\, \begin{aligned}&\exists\, x\in H^1(\intv, \R^n)\\ &\text{satisfying \eqref{eq:sys} a.e.\ in }\intv \end{aligned}\right\}.
\]
Throughout, we assume that \eqref{eq:sys} is minimal, i.e., $(A,B)$ is controllable and $(A,C)$ is observable. Observability is no restriction, as $\B$ depends only on the observable part of \eqref{eq:sys}. The \emph{lag} of $\B$, denoted $\lag(\B)$, is the smallest degree among all polynomial matrices $R$ such that, for all $(y,u)\in L^2(\intv, \R^p)\times L^2(\intv,\R^m)$, we have $(y,u)\in\B$ if and only if $R(\tfrac{\mathrm d}{\mathrm dt}) \col(y,u)=0$. By observability, $\lag(\B)$ equals the observability index of system~\eqref{eq:sys}, i.e., the smallest positive integer $k$ such that the Kalman observability matrix
\begin{equation*}
    \calO_k = \begin{bmatrix}
        C\\CA\\\vdots\\CA^{k-1}
    \end{bmatrix}\in\R^{kp\times n}
\end{equation*}
satisfies $\rank \calO_k = n$.

To obtain a data-based description of $\B$, the input used to generate the data must be sufficiently rich. This is measured by the Gramian of its jet. For $f\in H^{k-1}(\intv,\R^d)$ we define the Gramian
$$\Gamma_k(f) := \int_\intv \big(\Lambda_k(f)(t)\big)\big(\Lambda_k(f)(t)\big)^\top\diff t\in\R^{kd\times kd}.$$
For two functions $f\in H^{k-1}(\intv,\R^{d})$, $g\in H^{j-1}(\intv,\R^{e})$ the Gramian $\Gamma_{k,j}(f,g)$ is defined in a similar way with the stacked jet $\Lambda_{k,j}(f,g)$.
\begin{definition}[Persistency of excitation]\label{def:pe}
    For positive $L \in \mathbb{N}$, a signal $u:\intv\rightarrow \R^m$ is called persistently exciting of order $L$, if
    \begin{itemize}
        \item[(i)] $u\in H^{L-1}(\intv,\R^m)$,
        \item[(ii)] $\Gamma_L(u)$ is positive definite.
    \end{itemize}
\end{definition}

Using input--output trajectory data collected in a Gramian, the following version of the fundamental lemma for continuous-time LTI systems provides a data-based parameterization of the trajectories of $\B$.
\begin{lemma}[Fundamental lemma]
    \label{lem:fl}
    Let $L\in\N$ with $L\geq \lag(\B) +1$. Suppose that $(y,u)\in\B$ is given, where $u$ is persistently exciting of order $L+n$
as in \cref{def:pe}, and define $$\Gamma := \Gamma_{L,L}(y,u)\in\R^{L(p+m)\times L(p+m)}.$$ Then $(\bar y,\bar u)\in \B\cap( H^{L-1}(\intv,\R^p)\times H^{L-1}(\intv,\R^m))$ if and only if there exists $g\in L^2(\intv,\R^{L(p+m)})$ such that
    \begin{equation}
    \label{eq:FL}
        \Lambda_{L,L}(\bar y,\bar u) = \Gamma g.
    \end{equation}
    Furthermore, $\rank \Gamma=Lm + n$.
\end{lemma}

\paragraph{The fundamental lemma as a DAE}
Partition the data matrix~$\Gamma$ conformally with the jet $\Lambda_{L,L}(y,u)$ into block rows
\[
    \Gamma=\col\bigl(\Gamma_{y^{(0)}},\dots,\Gamma_{y^{(L-1)}},\Gamma_{u^{(0)}},\dots,\Gamma_{u^{(L-1)}}\bigr)
\]
with $\Gamma_{y^{(i)}}\in\R^{p\times L(p+m)}$ and $\Gamma_{u^{(i)}}\in\R^{m\times L(p+m)}$, $i=0,\dots,L-1$, and set
\begin{equation}
\label{eq:EandF}
\begin{aligned}
    \calE&:=\col\bigl(\Gamma_{y^{(0)}},\dots,\Gamma_{y^{(L-2)}},\Gamma_{u^{(0)}},\dots,\Gamma_{u^{(L-2)}}\bigr),\\
    \calF&:=\col\bigl(\Gamma_{y^{(1)}},\dots,\Gamma_{y^{(L-1)}},\Gamma_{u^{(1)}},\dots,\Gamma_{u^{(L-1)}}\bigr),
\end{aligned}
\end{equation}
so that $\calE,\calF\in\R^{(L-1)(p+m)\times L(p+m)}$.
Since consecutive block rows of~\eqref{eq:FL} are related by differentiation, \eqref{eq:FL} holds if and only if
$\bar y=\Gamma_{y^{(0)}}g$, $\bar u=\Gamma_{u^{(0)}}g$, and $g$ solves the differential-algebraic equation (DAE)
\begin{equation}
\label{eq:dae}
    \tfrac{\diff}{\diff t}(\calE g) = \calF g,
\end{equation}
i.e., $g\in L^2(\intv,\R^{L(p+m)})$ with $\calE g\in H^1(\intv,\R^{(L-1)(p+m)})$ and~\eqref{eq:dae} a.e.\ in $\intv$.

\paragraph{Problem statement}
\Cref{lem:fl} describes $\B$ only implicitly: the parameter $g$ is constrained by the rectangular, hence singular, DAE~\eqref{eq:dae}. Moreover, $g$ is not unique, since $\dim\ker\Gamma=Lp-n$.
We seek an explicit description of the full behavior $\B$
by a polynomial \emph{image matrix}
$M\in\R[s]^{(p+m)\times m}$ such that
$(\bar y,\bar u)\in\B$ if and only if
$$\col(\bar y,\bar u)=M\bigl(\tfrac{\diff}{\diff t}\bigr)h$$
for some latent variable $h\in L^2(\intv, \R^m)$. In contrast to $g$, the latent variable $h$ is unconstrained, and it is uniquely determined by $(\bar y,\bar u)$ if $M(\lambda)$ has full column rank for all $\lambda\in\C$.
In this case, $h$ is a flat output of~\eqref{eq:sys} in the sense of~\cite{fliessFlatnessDefectNonlinear1995}. 
Our objective is to compute such an $M$ directly from
$\Gamma$ using an orthogonal staircase reduction.
\section{Structure of the data pencil}
\label{sec:pencil}

The Gramian-based representation \eqref{eq:FL} leads to the
differential-algebraic equation \eqref{eq:dae}. We call its
associated pencil
\[
    s\calE-\calF
\]
the \emph{data pencil}. We first determine its Kronecker
structure. The main point is that, up to strict
equivalence, its nontrivial part consists of the state
equation augmented by the input differentiation chain.

For $k\geq1$, let
\begin{equation}
\label{eq:Toeplitz}
    \calT_k=
    \begin{bmatrix}
        D & 0 & \cdots & 0\\
        CB & D & \ddots & \vdots\\
        CAB & CB & \ddots & 0\\
        \vdots & \vdots & \ddots & D\\
        CA^{k-2}B & CA^{k-3}B & \cdots & CB & D
    \end{bmatrix}
    \in\R^{kp\times km}
\end{equation}
and
\[
    \calS_k=
    \begin{bmatrix}
        \calO_k & \calT_k\\
        0 & I_{km}
    \end{bmatrix}
    \in\R^{k(p+m)\times(n+km)}.
\]
Then every $(y,u)\in\B$ with $u\in H^{k-1}(\intv,\R^m)$ and corresponding state $x$ satisfies
$\Lambda_{k,k}(y,u)=\calS_k\Lambda_{1,k}(x,u)$.

\begin{lemma}
\label{lem:GammaFactor}
Let the assumptions of \cref{lem:fl} hold. 
Let $x$ be the state corresponding to the trajectory
$(y,u)$ and put $\Gamma_0:=\Gamma_{1,L}(x,u)$.  Then $\Gamma_0$ is positive
definite and
\begin{equation}
\label{eq:GammaFactor}
    \Gamma=\calS_L\Gamma_0\calS_L^\top.
\end{equation}
In particular, $\calS_L$ has full column rank and
$\Gamma_0\calS_L^\top$ has full row rank.
\end{lemma}

\begin{proof}
Positive definiteness of $\Gamma_0$ follows from
\cite[Prop.~21]{schmitzContinuoustimeFundamentalLemma2024}, and
\eqref{eq:GammaFactor} is equation~(50) in the proof of
\cite[Thm.~22]{schmitzContinuoustimeFundamentalLemma2024}.
Since $L\geq\lag(\B)+1$, observability implies that $\calO_L$ has full
column rank, and hence so does $\calS_L$.
\end{proof}

\subsection{Minimal bases and reduction to the augmented controllability pencil}
We recall only the facts on minimal bases that will be used below; see
\cite{forneyMinimalBases1975,vanDoorenKronecker1979}.  For a polynomial
matrix $P(s)\in\C[s]^{r\times q}$, its right rational nullspace is
$\ker_{\C(s)}P(s)$, and analogously for the left nullspace.  A polynomial
basis of such a rational subspace is called a \emph{minimal basis} if the
sum of its vector degrees is minimal.  The degrees of a minimal basis are
independent of the particular basis and are called the corresponding
\emph{minimal indices}.

We use the following standard characterization
\cite[Main Theorem~(2), p.~495]{forneyMinimalBases1975}.  Let
$N(s)=[n_1(s)\ \cdots\ n_k(s)]$ be polynomial and let $N_{\rm hc}$ be the
matrix of highest-degree column coefficients.  Then the
columns of $N$ form a minimal basis of the rational subspace they span if
and only if $N(\lambda)$ has full column rank for every $\lambda\in\C$ and
$N_{\rm hc}$ has full column rank.

For later reference, the right singular Kronecker block associated with a
minimal index $\varepsilon$ is
\[
    \calL_\varepsilon(s)=
    \begin{bmatrix}
        -1&s\\
          &-1&s\\
          &&\ddots&\ddots\\
          &&&-1&s
    \end{bmatrix}
    \in\C[s]^{\varepsilon\times(\varepsilon+1)}.
\]
Its right minimal index is $\varepsilon$; analogously,
$\calL_\eta(s)^\top$ represents a left minimal index $\eta$
\cite{vanDoorenKronecker1979}.

Define
\begin{equation}
\label{eq:PL}
    \calP_L(s)=
    \begin{bmatrix}
        sI_n-A & -B & 0 & \cdots & 0\\
        0 & sI_m & -I_m & & 0\\
        0 & 0 & sI_m & \ddots & \\
        \vdots & \vdots & & \ddots & -I_m
    \end{bmatrix}
    \in\R[s]^{(n+(L-1)m)\times(n+Lm)},
\end{equation}
with block columns corresponding to
$(x,u^{(0)},u^{(1)},\ldots,u^{(L-1)})$.
We call $\calP_L(s)$ the \emph{augmented controllability
pencil}. Its first block row contains the controllability
pencil
$\calR(s):=[\,sI_n-A\ \ -B\,]$.
The remaining block rows encode the input-jet relations
$su^{(j)}-u^{(j+1)}=0$, $j=0,\ldots,L-2$.

\begin{proposition}
\label{prop:pencil-equivalence}
Under the assumptions of \cref{lem:GammaFactor}, the data pencil $s\calE-\calF$ is
strictly equivalent to
\begin{equation}
\label{eq:pencil-equivalence}
    \begin{bmatrix}
        \calP_L(s)&0\\
        0&0
    \end{bmatrix},
\end{equation}
where the right zero block has $Lp-n$ columns and the lower one
has $(L-1)p-n$ rows.
\end{proposition}

\begin{proof}
Put $H:=\Gamma_0\calS_L^\top$, which has full row rank by
\cref{lem:GammaFactor}.  From the definitions of $\calE$ and $\calF$,
\[
s\calE-\calF=
\left(
s\begin{bmatrix}
    \calO_{L-1}&\calT_{L-1}&0\\
    0&I_{(L-1)m}&0
\end{bmatrix}
-
\begin{bmatrix}
    \calO_{L-1}A&\calO_{L-1}B&\calT_{L-1}\\
    0&0&I_{(L-1)m}
\end{bmatrix}
\right)H.
\]
The Toeplitz shift structure in \eqref{eq:Toeplitz} gives
\[
s\calE-\calF
    =
   \begin{bmatrix}
        I_{(L-1)p}&-\calT_{L-1}\\
        0&I_{(L-1)m}
    \end{bmatrix} \begin{bmatrix}
        \calO_{L-1}&0\\
        0&I_{(L-1)m}
    \end{bmatrix}
    \calP_L(s)H.
\]
Since $L-1\geq\lag(\B)$, the left factor multiplying $\calP_L$ has full
column rank.  Hence constant invertible row transformations reduce it to
$[\,I\ \ 0\,]^\top$.  Likewise, constant invertible column transformations
reduce the full-row-rank matrix $H$ to $[\,I\ \ 0\,]$.  This yields
\eqref{eq:pencil-equivalence}.
\end{proof}

The reduction separates the dynamics from two kinds of redundancy:
$\calP_L$ contains the state equation and the input jet, while the zero
columns describe nonuniqueness of the parameter $g$ and the zero rows
describe redundant output equations.

\subsection{Complete Kronecker structure}

Let $\nu_1,\ldots,\nu_m$ denote the controllability
indices of $(A,B)$, including zero indices when
necessary, so that $\sum_{j=1}^m\nu_j=n$.
Equivalently, they are the right minimal indices
of the controllability pencil $\calR(s)$;
see \cite{vanDoorenGeneralizedEigenstructure1981}.

\begin{theorem}
\label{thm:complete-Kronecker}
Under the assumptions of \cref{lem:GammaFactor}, the right minimal indices
of $s\calE-\calF$ are
\[
    \nu_1+L-1,\ldots,\nu_m+L-1
\]
together with $Lp-n$ zero indices.  Its left minimal indices are
$(L-1)p-n$ zero indices.  The pencil has no finite or infinite elementary
divisors.
\end{theorem}

\begin{proof}
By \cref{prop:pencil-equivalence}, it suffices to consider
the augmented controllability pencil $\calP_L$
from \eqref{eq:PL}.
Let
\[
    N_{\rm c}(s)=
    \begin{bmatrix}X(s)\\W(s)\end{bmatrix}
    \in\C[s]^{(n+m)\times m}
\]
be a minimal basis of $\ker_{\C(s)}\calR(s)$ with column degrees
$\nu_1,\ldots,\nu_m$. For each column $[x^\top\ w^\top]^\top$,
\[
    (sI_n-A)x(s)=Bw(s).
\]
If the column has degree $\nu$, comparison of the coefficient of
$s^{\nu+1}$ shows that the highest-degree coefficient of its state part
vanishes. Hence $\deg x<\nu$ and $\deg w=\nu$, where the degree of
the zero polynomial is understood as $-\infty$.

Now set
\[
    \widehat N_L(s)=
    \begin{bmatrix}
        X(s)\\
        W(s)\\
        sW(s)\\
        \vdots\\
        s^{L-1}W(s)
    \end{bmatrix}.
\]
Then $\calP_L(s)\widehat N_L(s)=0$.  Conversely, the lower block rows of
$\calP_L(s)z=0$ imply successively
$w_j=s^jw_0$, while the first block row gives
$(sI_n-A)x-Bw_0=0$.  Thus $\widehat N_L$ spans the complete right
rational nullspace of $\calP_L$.

Since $N_{\rm c}(\lambda)$ has full column rank for every
$\lambda$, so does $\widehat N_L(\lambda)$.
Moreover, its highest-column-degree coefficient matrix
is obtained from that of $N_{\rm c}$ by placing the
leading coefficients of $W$ in the last input block,
and therefore has full column rank.
Consequently $\widehat N_L$ is a minimal basis.  Its column degrees are
$\nu_j+L-1$.

The zero columns and rows in \eqref{eq:pencil-equivalence} contribute,
respectively, $Lp-n$ right and $(L-1)p-n$ left minimal indices equal to
zero.  Finally, the $m$ nontrivial right singular blocks use
$\sum_j(\nu_j+L-1)=n+(L-1)m$ rows.  Together with the $(L-1)p-n$ zero
left blocks this exhausts all rows, while the corresponding right blocks
together with the $Lp-n$ zero right blocks exhaust all columns.  Hence no
regular Kronecker block remains.
\end{proof}

Equivalently,
\[
s\calE-\calF
\sim
\operatorname{diag}\left(
    \calL_{\nu_1+L-1},\ldots,\calL_{\nu_m+L-1},
    \underbrace{\calL_0,\ldots,\calL_0}_{Lp-n},
    \underbrace{\calL_0^\top,\ldots,\calL_0^\top}_{(L-1)p-n}
\right),
\]
where $\calL_0$ denotes the empty $0\times1$ pencil.

\begin{remark}
\label{rem:common-kernels}
The zero minimal indices have a direct data interpretation. Since every block row of $\Gamma$ is a block row of $\calE$ or $\calF$,
\[
    \ker\Gamma
    =
    \ker
    \begin{bmatrix}
        \calE\\\calF
    \end{bmatrix}.
\]
This kernel has dimension $Lp-n$. 
By \cref{thm:complete-Kronecker},
$\dim\ker\col(
\calE^\top,\calF^\top)=(L-1)p-n$. 
Thus the zero right and left minimal indices can be extracted directly,
without a staircase iteration.
\end{remark}

\section{Orthogonal staircase reduction and image representation}
\label{sec:staircase}

The Kronecker form is useful for analysis but should not be computed
explicitly.  Staircase algorithms instead reveal the same structural
information by orthogonal transformations and rank-revealing
factorizations
\cite{vanDoorenKronecker1979,oaraImprovedAlgorithmComputation1997a}.
We exploit the special Kronecker structure of the data pencil
to simplify the staircase construction.

\subsection{The staircase of the data pencil}

A right staircase reduction reveals the right minimal indices
and possible infinite elementary divisors of a pencil by
orthogonal row and column transformations; see
\cite{vanDoorenKronecker1979,oaraImprovedAlgorithmComputation1997a}.

Let $s_i$ be the nullity of the active coefficient of $s$
at step $i$ and $r_i$ the rank of the compressed constant
block. Then
\[
    s_i-r_i=e_i,\qquad r_i-s_{i+1}=d_i,
\]
where $e_i$ is the multiplicity of the right minimal index
$i-1$ and $d_i$ counts infinite elementary divisors of
degree $i$; see \cite{vanDoorenKronecker1979}.

For the data pencil, \cref{thm:complete-Kronecker} gives
$d_i=0$, so $r_i=s_{i+1}$. Since there is no regular part
and all left minimal indices are zero, the terminal left
singular structure consists of $(L-1)p-n$ zero rows.

Thus, up to an orthogonal permutation of staircase blocks,
\begin{equation}
\label{eq:data-staircase}
\begin{split}
&P(sE-F)Q=\\
&\left[
\begin{array}{ccccc}
-A_{11}
 & sE_{12}-A_{12}
 & \cdots
 & sE_{1,\ell-1}-A_{1,\ell-1}
 & sE_{1\ell}-A_{1\ell}\\
0
 & -A_{22}
 & \cdots
 & sE_{2,\ell-1}-A_{2,\ell-1}
 & sE_{2\ell}-A_{2\ell}\\
\vdots & \ddots & \ddots & \vdots & \vdots\\
0 & \cdots & 0
 & -A_{\ell-1,\ell-1}
 & sE_{\ell-1,\ell}-A_{\ell-1,\ell}\\ \hline
0&\cdots&0&0&0
\end{array}
\right],
\end{split}
\end{equation}
where the $i$th block column has $s_i$ columns, the $i$th nonzero block
row has $s_{i+1}$ rows, $A_{ii}$ has full row rank, and
$E_{i,i+1}\in\R^{s_{i+1}\times s_{i+1}}$ is nonsingular.

Before deflation, the staircase block widths are completely
prescribed by the right minimal indices:
\[
    s_1=Lp-n+m,\qquad s_2=\cdots=s_L=m,
\]
and
\[
    s_{L+k}=\#\{j:\nu_j\geq k\},\qquad k\geq1.
\]
In particular, $\ell=L+\max_j\nu_j$.  Thus the first staircase step
contains the $Lp-n$ redundant directions of $g$; after they disappear,
all $m$ dynamic chains persist for another $L-1$ levels, and only then do
the chain lengths distinguish the controllability indices.

\subsection{A structure-exploiting orthogonal reduction}

The known zero minimal indices can be removed before the recursive
staircase is computed. This reduces the dimensions of the matrices
in the subsequent staircase steps.

Let $Z_0$ and $Y_0$ have orthonormal columns spanning
\[
    \ker\begin{bmatrix}E\\F\end{bmatrix},
    \qquad
    \ker\begin{bmatrix}E^\top\\F^\top\end{bmatrix},
\]
respectively, and complete them to orthogonal matrices
$Z=[\,Z_0\ Z_1\,]$ and $Y=[\,Y_1\ Y_0\,]$.  By
\cref{rem:common-kernels}, $Z_0$ has $Lp-n$ columns and $Y_0$ has
$(L-1)p-n$ columns.  Hence
\begin{equation}
\label{eq:deflated-pencil}
Y^\top(sE-F)Z=
\begin{bmatrix}
    0&sE_1-F_1\\
    0&0
\end{bmatrix},
\qquad
E_1=Y_1^\top EZ_1,\quad F_1=Y_1^\top FZ_1.
\end{equation}
The reduced pencil has size
$(n+(L-1)m)\times(n+Lm)$ and consists exclusively of the $m$ right
singular chains with minimal indices $\nu_j+L-1$.

The staircase in \eqref{eq:data-staircase} can now
be computed recursively. Suppose the
active pair at step $i$ is $(E_i,F_i)$.  Compute an orthonormal basis
$Z_i^0$ of $\ker E_i$, set $s_i=\dim\ker E_i$, and complete
$Z_i=[\,Z_i^0\ Z_i^1\,]$ orthogonally.  Next form
$C_i=F_iZ_i^0$, set $r_i=\rank C_i$, and choose an orthogonal
$Y_i=[\,Y_i^{\rm a}\ Y_i^{\rm b}\,]$ such that
\[
    Y_i^\top C_i=
    \begin{bmatrix}
        A_{ii}\\0
    \end{bmatrix},
    \qquad \rank A_{ii}=r_i.
\]
Then
\begin{equation}
\label{eq:staircase-step}
Y_i^\top(sE_i-F_i)Z_i=
\begin{bmatrix}
    -A_{ii}
    &Y_i^{{\rm a}\top}(sE_i-F_i)Z_i^1\\
    0&sE_{i+1}-F_{i+1}
\end{bmatrix},
\end{equation}
where
$E_{i+1}=Y_i^{{\rm b}\top}E_iZ_i^1$ and
$F_{i+1}=Y_i^{{\rm b}\top}F_iZ_i^1$.
The process stops when the active pencil has no rows;
the number of its remaining columns is recorded as the
final staircase width.

For the present pencil no infinite elementary divisors are possible, so
$r_i=s_{i+1}$ throughout. For the reduced pencil obtained after the deflation
\eqref{eq:deflated-pencil}, the staircase widths satisfy
$s_1=\cdots=s_L=m$. Consequently, the first $L-1$ compressed matrices
$C_i$ have full column rank $m$. Their ranks are thus known
in exact arithmetic, while rank loss from level $L$ onward
determines the controllability indices.

\paragraph{Algorithm 1: Structure-exploiting staircase reduction}
Given $E=\calE$ and $F=\calF$:
\begin{enumerate}
    \item Compute orthonormal bases $Z_0$ of
        $\ker\col(E, F)$ and $Y_0$ of
    $\ker\col(E^\top, F^\top)$.
    Complete them orthogonally and form the reduced pair $(E_1,F_1)$ in
    \eqref{eq:deflated-pencil}.

    \item For the current pair $(E_i,F_i)$, compute an orthonormal basis
    $Z_i^0$ of $\ker E_i$ and an orthogonal completion
    $Z_i=[\,Z_i^0\ Z_i^1\,]$.

    \item Form $C_i=F_iZ_i^0$.  Compute a rank-revealing orthogonal
    factorization of $C_i$ and an orthogonal row transformation
    $Y_i=[\,Y_i^{\rm a}\ Y_i^{\rm b}\,]$ for which
    $Y_i^\top C_i=[\,A_{ii}^\top\ 0\,]^\top$.

    \item Form the next active pair from the lower-right block in
    \eqref{eq:staircase-step}, record $s_i=\dim\ker E_i$ and
    $r_i=\rank C_i$, and continue until no rows remain,
recording the final number of columns.

        \item Accumulate all orthogonal transformations and check
    $r_i=s_{i+1}$, as required in exact arithmetic. Violations suggest
    inaccurate numerical rank estimates or a departure from the assumed
    pencil structure.
\end{enumerate}

All transformations in this procedure are orthogonal.  In particular,
the full-rank staircase blocks are not normalized to identity matrices.

\subsection{From the staircase to a polynomial kernel basis}

The parameter $g$ in \eqref{eq:FL} is constrained by the DAE \eqref{eq:dae}. In this subsection, we use the staircase to compute a minimal polynomial basis of the right nullspace of $s\calE-\calF$. This basis turns the implicit constraint into an explicit parameterization of all solutions $g$ by free latent variables, which is the basis for the image representation in \Cref{sec:image}.

After the deflation \eqref{eq:deflated-pencil}, let
\[
    P_1(sE_1-F_1)Q_1=S(s){}=
    \begin{bmatrix}
    -A_{11} & sE_{12}-A_{12} & \cdots & sE_{1\ell}-A_{1\ell}\\
    0 & -A_{22} & \cdots & sE_{2\ell}-A_{2\ell}\\
    \vdots & \ddots & \ddots & \vdots\\
    0 & \cdots & -A_{\ell-1,\ell-1} & sE_{\ell-1,\ell}-A_{\ell-1,\ell}
    \end{bmatrix}
\]
be the pure right-singular staircase obtained above, i.e., \eqref{eq:data-staircase} without zero rows.  Write its block
column sizes as $s_1,\ldots,s_\ell$ and set $s_{\ell+1}=0$ and
$e_i=s_i-s_{i+1}$.  Thus $\sum_i e_i=m$, and $e_i$ is the number of
dynamic right minimal indices equal to $i-1$.

A polynomial kernel vector $\col(n_1,\ldots,n_\ell)$ of $S(s)$ can be
computed by back substitution. Its last block $n_\ell$ is free, and the
$i$th block row requires
$A_{ii}n_i(s)=\sum_{j>i}(sE_{ij}-A_{ij})n_j(s)$. Since $A_{ii}$ has full row
rank, $n_i$ is determined by $n_{i+1},\ldots,n_\ell$ up to an element
of $\ker A_{ii}$, which has dimension $e_i$. We therefore need a right
inverse and a kernel basis of $A_{ii}$.
For $i<\ell$, choose an orthogonal matrix
$V_i=[\,V_i^{\rm r}\ V_i^0\,]$ such that
\[
    A_{ii}V_i=\begin{bmatrix}R_i&0\end{bmatrix},
\]
where $R_i$ is square and nonsingular.  Put
$G_i=V_i^{\rm r}R_i^{-1}$ and $K_i=V_i^0$, so that $A_{ii}G_i=I$ and $A_{ii}K_i=0$; numerically, multiplication
with $R_i^{-1}$ is implemented by solving with $R_i$.

Let $h=\col(h_1,\ldots,h_\ell)\in\R^m$, with
$h_i\in\R^{e_i}$, and let $\Pi_i h=h_i$.  Define polynomial matrices
$N_i(s)$ backwards by
\begin{equation}
\label{eq:minimal-basis-recursion}
    N_\ell(s)=\Pi_\ell,\qquad
    N_i(s)=K_i\Pi_i+
    G_i\sum_{j=i+1}^{\ell}(sE_{ij}-A_{ij})N_j(s).
\end{equation}
Then
\[
    N_{\rm st}(s)=
    \col\bigl(N_1(s),\ldots,N_\ell(s)\bigr)
\]
spans the right rational nullspace of $S(s)$.  The staircase
reconstruction results of
\cite{noferiniVanDoorenMinimalBases2023} imply that its columns form a
minimal basis, with $e_i$ columns of degree $i-1$.

Transforming back gives
\begin{equation}
\label{eq:dynamic-minimal-basis}
    N(s):=Z_1Q_1N_{\rm st}(s)
    \in\R[s]^{L(p+m)\times m}.
\end{equation}
Consequently,
\begin{equation}
\label{eq:full-kernel-basis}
    \ker_{\R(s)}(s\calE-\calF)
    =
    \operatorname{im}_{\R(s)}
    \begin{bmatrix}
        Z_0&N(s)
    \end{bmatrix}.
\end{equation}
Since $Z_0$ and $Z_1Q_1$ have orthogonal ranges,
the pointwise and highest-column-degree rank conditions
for $N_{\rm st}$ carry over to $[\,Z_0\ N(s)\,]$.
Thus the latter is a minimal basis of the full right
nullspace, with $Lp-n$ columns of degree zero and,
after a suitable ordering, the remaining column degrees
$\nu_1+L-1,\ldots,\nu_m+L-1$.

The polynomial basis in \eqref{eq:full-kernel-basis}
also parameterizes all weak solutions of \eqref{eq:dae}. Indeed, the reduced staircase pencil $S(s)$ has full
row rank for every $s\in\C$. Hence it admits a polynomial right inverse $R(s)$.
Since $N_{\rm st}(s)$ is a minimal basis, its columns also form
a basis of the polynomial kernel of $S(s)$. Consequently,
$U(s)=[\,N_{\rm st}(s)\ R(s)\,]$ is unimodular. Write $U^{-1}=\col(U_1,U_2)$. Then $SU=[\,0\ I\,]$ gives $U_2=S$, and $U^{-1}U=I$ gives $U_1N_{\rm st}=I$.
Consequently, every distributional solution of
$S(\tfrac{\diff}{\diff t})z=0$ has a unique representation
\[
 z=N_{\rm st}\!\left(\tfrac{\diff}{\diff t}\right)h
\]
with distributional components of $h$, namely $h=U_1(\tfrac{\diff}{\diff t})z$.

We investigate the regularity of $h$ given $z\in L^2(\intv)$. Write
\[
 N_{\rm st}(s)
 =
 \begin{bmatrix}n_1(s)&\cdots&n_m(s)\end{bmatrix},
 \qquad
 n_j(s)=\sum_{k=0}^{\varepsilon_j}n_{jk}s^k,
\]
and let
\[
 A_{\rm hc}
 =
 \begin{bmatrix}
 n_{1,\varepsilon_1}&\cdots&n_{m,\varepsilon_m}
 \end{bmatrix}.
\]
Since $N_{\rm st}$ is a minimal basis, $A_{\rm hc}$
has full column rank. Choose a constant left inverse
$L_{\rm hc}$ of $A_{\rm hc}$.

For $z=N_{\rm st}(\tfrac{\diff}{\diff t})h$, we obtain,
in the sense of distributions,
\[
 \begin{bmatrix}
 h_1^{(\varepsilon_1)}\\
 \vdots\\
 h_m^{(\varepsilon_m)}
 \end{bmatrix}
 =
 L_{\rm hc}z
 -
 L_{\rm hc}
 \sum_{j=1}^m\sum_{k=0}^{\varepsilon_j-1}
 n_{jk}h_j^{(k)}.
\]
Since all column degrees satisfy
$\varepsilon_j\geq L-1\geq1$, let $v$ collect
$h_j,\ldots,h_j^{(\varepsilon_j-1)}$ for all $j$.
The preceding identity yields a constant-coefficient
system
\[
 \dot v=A_0v+B_0z
\]
in the sense of distributions. Every distributional
solution is given by the variation-of-constants formula.
Since $\intv$ is bounded and $z\in L^2(\intv)$,
we obtain $v\in H^1(\intv)$ and hence
$h_j\in H^{\varepsilon_j}(\intv)$ for all $j$.
The converse follows directly by differentiation.

Transforming back and restoring the freely varying component in
$\ker\Gamma$, every $L^2$-solution of \eqref{eq:dae}
therefore has the representation
\[
 g=Z_0h_0+
 N\!\left(\tfrac{\diff}{\diff t}\right)h,
\]
where $h_0\in L^2(\intv,\R^{Lp-n})$ and, after ordering
according to the controllability indices,
$h_j\in H^{\nu_j+L-1}(\intv)$.

\subsection{Data-driven image representation}\label{sec:image}

Set
\[
    G:=
    \begin{bmatrix}
        \Gamma_{y^{(0)}}\\
        \Gamma_{u^{(0)}}
    \end{bmatrix}
\]
and define
\begin{equation}
\label{eq:image-matrix}
    M(s):=GN(s)=
    \begin{bmatrix}
        M_y(s)\\
        M_u(s)
    \end{bmatrix}
    \in\R[s]^{(p+m)\times m}.
\end{equation}
Since $\operatorname{im}Z_0=\ker\Gamma$, the arbitrary term $Z_0h_0$ is
invisible in the manifest variables.  Moreover,
$(s\calE-\calF)N(s)=0$ and the shift structure of $\calE,\calF$ imply
\begin{equation}
\label{eq:jet-factorization}
    \Gamma N(s)=
    \begin{bmatrix}
        M_y(s)\\
        sM_y(s)\\
        \vdots\\
        s^{L-1}M_y(s)\\
        M_u(s)\\
        sM_u(s)\\
        \vdots\\
        s^{L-1}M_u(s)
    \end{bmatrix}.
\end{equation}

\begin{theorem}[Data-driven image representation]
\label{thm:data-image-representation}
After a suitable ordering of the columns, the column degrees of $M$ are
the controllability indices $\nu_1,\ldots,\nu_m$, and
$\rank M(\lambda)=m$ for every $\lambda\in\C$.  Furthermore,
\begin{equation}
\label{eq:behavior-image}
\begin{split}
&\B\cap
\bigl(
H^{L-1}(\intv,\R^p)\times H^{L-1}(\intv,\R^m)
\bigr)
\\
&\quad=
\left\{
(\bar y,\bar u)\ \middle|\
\begin{bmatrix}\bar y\\\bar u\end{bmatrix}
=
M\!\left(\tfrac{\diff}{\diff t}\right)h,\quad
h_j\in H^{\nu_j+L-1}(\intv,\R),\
j=1,\ldots,m
\right\}.
\end{split}
\end{equation}
\end{theorem}

\begin{proof}
Let $n_j(s)$ be a column of $N$ of degree
$\varepsilon_j=\nu_j+L-1$.  The matrix
$[\,Z_0\ N(s)\,]$ is a minimal basis of the full right nullspace, so its
highest-column-degree coefficient matrix has full column rank.  Hence the
leading coefficient of $n_j$ does not belong to
$\operatorname{im}Z_0=\ker\Gamma$, and therefore
$\deg(\Gamma n_j)=\varepsilon_j$.  By \eqref{eq:jet-factorization},
$\deg(\Gamma n_j)=\deg M_j+L-1$, giving $\deg M_j=\nu_j$.

If $M(\lambda)c=0$, then \eqref{eq:jet-factorization} gives
$\Gamma N(\lambda)c=0$, and hence
$N(\lambda)c\in\operatorname{im}Z_0$.  Since
$[\,Z_0\ N(\lambda)\,]$ has full column rank, $c=0$.
Since the maximal minors of $M$ have no common complex
zero, $M$ admits a polynomial left inverse. Thus $h$
is a differential function of the manifest variables.
Conversely, the identity
\[
 \Lambda_{\lag(\B)}(y)
 =
 \calO_{\lag(\B)}x+
 \calT_{\lag(\B)}\Lambda_{\lag(\B)}(u)
\]
and the full column rank of $\calO_{\lag(\B)}$
show that $x$ is a differential function of $(y,u)$.
Since $(y,u)$ is generated by $M(\tfrac{\diff}{\diff t})h$,
both $x$ and $u$ are differential functions of $h$.
Hence $h$ is a flat output of \eqref{eq:sys}.

Finally, if
$h_j\in H^{\nu_j+L-1}(\intv)$, then
$g=N(\tfrac{\diff}{\diff t})h$ satisfies \eqref{eq:dae}, and
\eqref{eq:jet-factorization} gives
\[
    \Gamma g=\Lambda_{L,L}(\bar y,\bar u),
    \qquad
    \begin{bmatrix}\bar y\\\bar u\end{bmatrix}
    =
    M\!\left(\tfrac{\diff}{\diff t}\right)h.
\]
Hence $(\bar y,\bar u)\in\B$ by \cref{lem:fl}.

Conversely, any trajectory in the left-hand side of
\eqref{eq:behavior-image} admits, by \cref{lem:fl}, a parameter $g$
satisfying \eqref{eq:dae}. The weak-solution parameterization established above gives
$g=Z_0h_0+N(\tfrac{\diff}{\diff t})h$, with the stated Sobolev
regularity. Since $\Gamma Z_0=0$, the first term disappears after
multiplication by $\Gamma$. This proves \eqref{eq:behavior-image}.
\end{proof}

\paragraph{Algorithm 2: Data-driven construction of the image matrix}
Starting from the data Gramian $\Gamma$:
\begin{enumerate}
    \item Form $\calE,\calF$ from \eqref{eq:EandF} and
    $G=[\,\Gamma_{y^{(0)}}^\top\ \Gamma_{u^{(0)}}^\top\,]^\top$.

    \item Apply Algorithm~1, first extracting $Z_0$ and $Y_0$ and then
    computing the orthogonal staircase of the reduced pencil.

    \item Reconstruct $N_{\rm st}(s)$ by
    \eqref{eq:minimal-basis-recursion} and transform back according to
    \eqref{eq:dynamic-minimal-basis}.

    \item Return
    \[
        M(s)=GN(s).
    \]
\end{enumerate}
The resulting matrix has $m$ columns, with column degrees
equal to the controllability indices. The higher regularity
in \eqref{eq:behavior-image} reflects the restriction to
$H^{L-1}$ in the fundamental lemma. The same image matrix
also represents the full behavior, as shown next.

\begin{corollary}[Image representation of the full behavior]
\label{cor:full-image}
Let $M$ be given by \eqref{eq:image-matrix}. Then
\[
\B=
\left\{
(\bar y,\bar u)\in
L^2(\intv,\R^p)\times L^2(\intv,\R^m)
\ \middle|\
\begin{bmatrix}\bar y\\\bar u\end{bmatrix}
=
M\!\left(\tfrac{\diff}{\diff t}\right)h,\quad
h_j\in H^{\nu_j}(\intv,\R),\
j=1,\ldots,m
\right\}.
\]
The latent variable $h$ is uniquely determined by
$(\bar y,\bar u)$.
\end{corollary}

\begin{proof}
Let $N_{\rm hc}$ and $M_{\rm hc}$ denote the matrices
formed by the highest-degree coefficient of each column
of $N$ and $M$, respectively. By
\eqref{eq:jet-factorization},
\[
\Gamma N_{\rm hc}=J M_{\rm hc},
\]
where
\[
J=\operatorname{diag}
\bigl(\mathbf e_L\otimes I_p,\mathbf e_L\otimes I_m\bigr)
\in\R^{L(p+m)\times(p+m)}
\]
and $\mathbf e_L$ is the last canonical basis vector
of $\R^L$.
Since $[\,Z_0\ N_{\rm hc}\,]$ has full column rank and
$\ker\Gamma=\operatorname{im}Z_0$, the matrix
$\Gamma N_{\rm hc}$ has full column rank.
Hence $M_{\rm hc}$ has full column rank.

Suppose that $h$ is distributional and
$w=M(\tfrac{\diff}{\diff t})h\in L^2(\intv)$.
A constant left inverse of $M_{\rm hc}$ allows us
to express the highest derivatives
$h_j^{(\nu_j)}$ in terms of $w$ and lower-order
derivatives of $h$. The lower-order derivatives
corresponding to $\nu_j>0$ therefore satisfy a
constant-coefficient first-order system with
$L^2$ forcing. The components with $\nu_j=0$
are recovered algebraically. The
variation-of-constants formula yields
$h_j\in H^{\nu_j}(\intv)$ for every $j$.

By the proof of \cref{thm:data-image-representation}, $M$ admits
a polynomial left inverse $P_{\rm l}(s)$.
Furthermore, $\B$ is closed in $L^2$: by variation
of constants, the output is a bounded linear
function of the input plus an element of a
finite-dimensional subspace of $L^2$.
Smooth trajectories are dense in $\B$, since
the input can be approximated by smooth functions
while keeping the initial state fixed.

Now let $w\in\B$ and choose smooth trajectories
$w_k\in\B$ with $w_k\to w$ in $L^2$.
By \cref{thm:data-image-representation}, there
exist smooth $h_k$ such that
\[
w_k=M\!\left(\tfrac{\diff}{\diff t}\right)h_k,
\qquad
h_k=P_{\rm l}\!\left(\tfrac{\diff}{\diff t}\right)w_k.
\]
Consider the test function space $\mathcal D$ consisting of infinitely differentiable $\varphi:\intv\to\R^{p+m}$ with compact support, and define $M^*(s):=M(-s)^\top$ and $P_{\rm l}^*(s):=P_{\rm l}(-s)^\top$. These are the formal adjoints of $M$ and $P_{\rm l}$, that is, integration by parts moves $M(\tfrac{\diff}{\diff t})$ and $P_{\rm l}(\tfrac{\diff}{\diff t})$ onto the test function as $M^*(\tfrac{\diff}{\diff t})$ and $P_{\rm l}^*(\tfrac{\diff}{\diff t})$. This yields
\[
    \int_\intv w_k(t)^\top \varphi(t) \, \mathrm dt = \int_\intv h_k(t)^\top (M^*(\tfrac{\mathrm d}{\mathrm dt}) \varphi)(t) \, \mathrm dt = \int_\intv w_k(t)^\top \big(P_{\rm l}^*(\tfrac{\mathrm d}{\mathrm dt}) M^*(\tfrac{\mathrm d}{\mathrm dt}) \varphi\big)(t)\,\mathrm dt
\]
for all $\varphi\in \mathcal D$. Passing to the limit $k\rightarrow \infty$ shows
\[
\int_\intv w(t)^\top \varphi(t) \, \mathrm dt =  \int_\intv w(t)^\top \big(P_{\rm l}^*(\tfrac{\mathrm d}{\mathrm dt}) M^*(\tfrac{\mathrm d}{\mathrm dt}) \varphi\big)(t)\,\mathrm dt
\]
for all $\varphi\in \mathcal D$. Since the right-hand side does not involve derivatives of $w$, this identity is meaningful for $w\in L^2$, and it states that the distribution
\[
h:=P_{\rm l}\!\left(\tfrac{\diff}{\diff t}\right)w
\]
satisfies $w=M(\tfrac{\diff}{\diff t})h$.
By the regularity argument above,
$h_j\in H^{\nu_j}(\intv)$.

Conversely, let $h_j\in H^{\nu_j}(\intv)$ and
approximate each $h_j$ by smooth functions in
this Sobolev space. By
\cref{thm:data-image-representation}, the resulting
smooth trajectories belong to $\B$. They converge
in $L^2$ to $M(\tfrac{\diff}{\diff t})h$.
Since $\B$ is closed, the limit also belongs
to $\B$. Uniqueness follows from
$P_{\rm l}(s)M(s)=I_m$.
\end{proof}

\section{Numerical aspects and validation}
\label{sec:numerics}

The construction of \cref{sec:staircase} relies on numerical rank decisions
for matrices formed from blocks of $\Gamma$. We first describe a scaling
designed to improve the conditioning of $\Gamma$ while preserving the
structure exploited by Algorithm~2 (\cref{subsec:preconditioning}),
and then illustrate the
complete procedure on a chain of mass-spring-damper systems
(\cref{subsec:example}).

\subsection{Preconditioning the Gramian}
\label{subsec:preconditioning}

The numerical conditioning of derivative Gramians depends strongly on the
relative scaling of derivatives of different orders.  This becomes particularly relevant when $L$ is moderately large, as may be
required by the lag of the behavior.

Since the Gramian $\Gamma$ is in general singular, we measure its conditioning
on its range.  For a nonzero positive semidefinite matrix $G$, let
\[
    \kappa_+(G)
    :=
    \frac{\lambda_{\max}(G)}{\lambda_{\min}^+(G)},
\]
where $\lambda_{\min}^+(G)$ denotes its smallest positive eigenvalue.  Under
the assumptions of \cref{lem:GammaFactor}, $\rank\Gamma=n+Lm$, so that
$\kappa_+(\Gamma)$ measures the separation of the nonzero singular values
relevant for the numerical rank decisions in the staircase algorithm.

The following observation describes the dependence of the Gramian on the time
scale.

\begin{lemma}
\label{lem:gramian-scaling}
Let $f\in H^{L-1}(\intv,\R^d)$, let $c>0$, and define
$f_c(t):=f(ct)$ for $t\in c^{-1}\intv$.  With
\[
    D_c
    :=
    \operatorname{diag}
    \bigl(
        c^{-1/2},c^{1/2},\ldots,c^{L-3/2}
    \bigr),
\]
one has
\begin{equation}
\label{eq:gramian-scaling}
    \Gamma_L(f_c)
    =
    (D_c\otimes I_d)\,
    \Gamma_L(f)\,
    (D_c\otimes I_d).
\end{equation}
\end{lemma}

\begin{proof}
For $k=0,\ldots,L-1$,
$f_c^{(k)}(t)=c^k f^{(k)}(ct)$.  Hence the $(k,j)$ block of
$\Gamma_L(f_c)$ is
\[
    \int_{c^{-1}\intv}
        c^{k+j}f^{(k)}(ct)f^{(j)}(ct)^\top\,\diff t
    =
    c^{k+j-1}
    \int_{\intv}
        f^{(k)}(t)f^{(j)}(t)^\top\,\diff t,
\]
which is precisely \eqref{eq:gramian-scaling}.
\end{proof}

Thus a change of time scale by a factor $c$ changes the relative scaling of
the zeroth and $(L-1)$st derivatives by a factor $c^{L-1}$, and that of the
corresponding diagonal Gramian blocks by $c^{2(L-1)}$.  For large $L$, even a
moderate mismatch in the time scale may therefore cause severe numerical
scaling problems.

Motivated by \cref{lem:gramian-scaling}, we consider preconditioners of the
form
\begin{equation}
\label{eq:preconditioner}
    P(c,T_{\mathrm y},T_{\mathrm u})
    :=
    \operatorname{diag}
    \bigl(
        D_c\otimes T_{\mathrm y},
        D_c\otimes T_{\mathrm u}
    \bigr),
\end{equation}
where $c>0$ and
$T_{\mathrm y}\in\R^{p\times p}$ and
$T_{\mathrm u}\in\R^{m\times m}$ are invertible.  The parameter $c$ balances
the different derivative orders, while $T_{\mathrm y}$ and
$T_{\mathrm u}$ balance the output and input channels, respectively.

An important feature of \eqref{eq:preconditioner} is that Algorithm~2 can be
applied directly to the preconditioned Gramian and the image representation
of the original behavior can subsequently be recovered explicitly.

\begin{lemma}
\label{lem:imrep-recovery}
Let
\[
    \widetilde\Gamma
    =
    P\Gamma P^\top,
    \qquad
    P=P(c,T_{\mathrm y},T_{\mathrm u}),
\]
and let $\widetilde M$ be the polynomial image matrix obtained by applying
Algorithm~2 to $\widetilde\Gamma$.  Then
\begin{equation}
\label{eq:imrep-recovery}
    M(s)
    :=
    \operatorname{diag}(T_{\mathrm y},T_{\mathrm u})^{-1}
    \widetilde M(cs)
\end{equation}
is an image representation of the original behavior $\B$.
\end{lemma}

\begin{proof}
Let
\[
    T:=\operatorname{diag}(T_{\mathrm y},T_{\mathrm u})
\]
and consider the transformed manifest signal
\[
    \begin{bmatrix}
        \widetilde y(t)\\
        \widetilde u(t)
    \end{bmatrix}
    =
    T
    \begin{bmatrix}
        y(ct)\\
        u(ct)
    \end{bmatrix},
    \qquad
    t\in c^{-1}\intv.
\]
This defines a bijection between $\B$ and the transformed behavior
$\widetilde\B$ on $c^{-1}\intv$.  Moreover,
\[
    \Lambda_{L,L}(\widetilde y,\widetilde u)(t)
    =
    \sqrt c\,
    P\,
    \Lambda_{L,L}(y,u)(ct).
\]
A change of variables therefore gives
\[
    \Gamma_{L,L}(\widetilde y,\widetilde u)
    =
    P\Gamma P^\top
    =
    \widetilde\Gamma.
\]
Consequently, Algorithm~2 applied to $\widetilde\Gamma$ yields
\[
    \begin{bmatrix}
        \widetilde y\\
        \widetilde u
    \end{bmatrix}
    =
    \widetilde M\!\left(\tfrac{\diff}{\diff t}\right)\widetilde h
\]
for the transformed behavior.

Given a free signal $h$ on $\intv$, set
$\widetilde h(t)=h(ct)$.  Evaluating the preceding representation at
$t/c$ gives
\[
    \begin{bmatrix}
        y(t)\\
        u(t)
    \end{bmatrix}
    =
    T^{-1}
    \widetilde M\!\left(
        c\tfrac{\diff}{\diff t}
    \right)h(t).
\]
Hence \eqref{eq:imrep-recovery} is an image representation of $\B$.
\end{proof}

It remains to choose the scaling parameters.  A natural objective is
\[
    \min_{c,T_{\mathrm y},T_{\mathrm u}}
    \kappa_+
    \bigl(
        P(c,T_{\mathrm y},T_{\mathrm u})
        \Gamma
        P(c,T_{\mathrm y},T_{\mathrm u})^\top
    \bigr).
\]
This problem is nonconvex, and we therefore use a simple alternating
heuristic.

For $k=0,\ldots,L-1$, the diagonal blocks of $\Gamma$ associated with
$y^{(k)}$ and $u^{(k)}$ are the Gramians $\Gamma_1(y^{(k)})\in\R^{p\times p}$
and $\Gamma_1(u^{(k)})\in\R^{m\times m}$.

\paragraph{Algorithm 3: preconditioner search}
Given $\Gamma$ and its diagonal derivative blocks:
\begin{enumerate}
    \item Initialize
    $c=1$, $T_{\mathrm y}=I_p$, $T_{\mathrm u}=I_m$, and retain the
    corresponding preconditioner as the current best choice.

    \item Repeat until convergence or a prescribed maximum number of
    iterations is reached:
    \begin{itemize}
        \item[(a)] With $T_{\mathrm y}$ and $T_{\mathrm u}$ fixed, perform a
        one-dimensional grid search in $c>0$ to reduce
        $\kappa_+(P\Gamma P^\top)$.

        \item[(b)] For the resulting $c$, form
        \[
            R_{\mathrm y}
            =
            \sum_{k=0}^{L-1}
                c^{2k-1}\Gamma_1(y^{(k)}),
            \qquad
            R_{\mathrm u}
            =
            \sum_{k=0}^{L-1}
                c^{2k-1}\Gamma_1(u^{(k)}).
        \]

        \item[(c)] If $R_{\mathrm y}$ and $R_{\mathrm u}$ are positive
        definite, set
        \[
            T_{\mathrm y}=R_{\mathrm y}^{-1/2},
            \qquad
            T_{\mathrm u}=R_{\mathrm u}^{-1/2}.
        \]

        \item[(d)] Retain the current preconditioner whenever it improves
        $\kappa_+(P\Gamma P^\top)$.
    \end{itemize}

    \item Return the best preconditioner found.
\end{enumerate}

The channel transformations in step~(c) whiten the aggregate derivative
energies for the current time scale, whereas step~(a) balances the relative
weight of the different derivative orders.  No claim of global optimality is
made; the purpose of the procedure is to obtain a well-scaled Gramian for the
subsequent numerical rank decisions.

\subsection{Numerical example}
\label{subsec:example}
\begin{figure}
    \centering
\scalebox{0.9}{
\begin{tikzpicture}
    % wall
    \draw (0,1) -- (0,-1);
    \fill[pattern=north west lines] (-0.2,-1) rectangle (0,1);
    % msd_1
    \draw (0,0.5) to[spring, l=$k_1$] (3,0.5);
    \draw (0,-0.5) to[damper, l_=$d_1$,] (3,-0.5);
    \node[mass, draw, fill=gray!30, minimum width=1cm, minimum height=1.5cm] (M1) at (3.5,0) {$m_1$};
    \draw[thick, -latex] (2,0) -- (3,0) node[midway, above] {$u_1$};
    % msd_2
    \draw (4,0.5) to[spring, l=$k_2$] (7,0.5);
    \draw (4,-0.5) to[damper, l_=$d_2$,] (7,-0.5);
    \node[mass, draw, fill=gray!30, minimum width=1cm, minimum height=1.5cm] (M1) at (7.5,0) {$m_2$};
    \draw[thick, -latex] (6,0) -- (7,0) node[midway, above] {$u_2$};
    % dots
    \draw[dotted, thick] (8,0.5) -- (8.5,0.5);
    \draw[dotted, thick] (8,-0.5) -- (8.5,-0.5);
    % msd_n
    \draw (8.5,0.5) to[spring, l=$k_N$] (11.5,0.5);
    \draw (8.5,-0.5) to[damper, l_=$d_N$,] (11.5,-0.5);
    \node[mass, draw, fill=gray!30, minimum width=1cm, minimum height=1.5cm] (Mn) at (12,0) {$m_N$};
    \draw[thick, -latex] (10.5,0) -- (11.5,0) node[midway, above] {$u_N$};
\end{tikzpicture}
}
\caption{Chain of coupled mass-spring-damper systems.}
    \label{fig:msd_chain}
\end{figure}
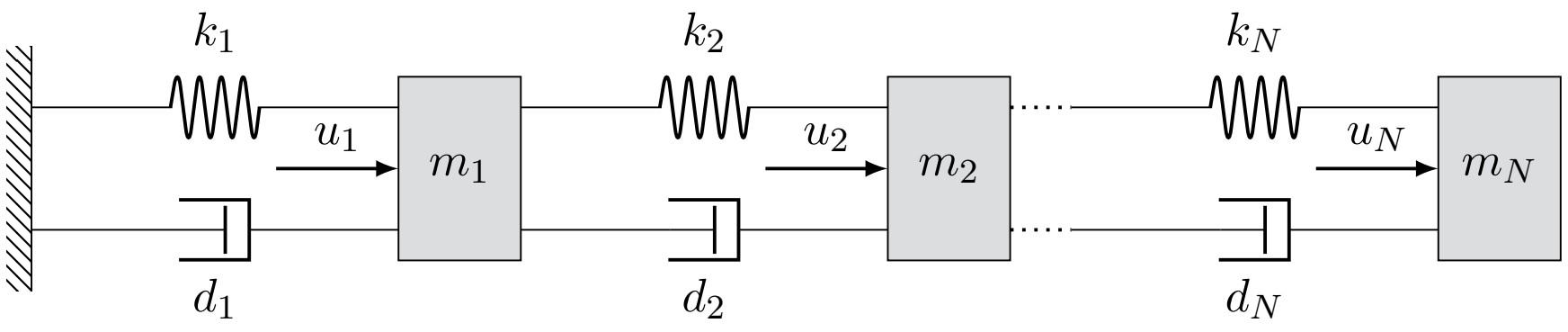
We demonstrate Algorithm~2 on a chain of $N$ coupled mass-spring-damper systems, see \cref{fig:msd_chain}. The state vector collects the displacement $q_i$ and the velocity $\dot q_i$ of each mass $m_i$, i.e., $x = [\, q_1,\dots, q_N,\dot q_1,\dots,\dot q_N\,]^\top$. Adjacent masses are coupled by spring-damper pairs with stiffness $k_i$ and damping $d_i$. The leftmost mass is attached to a wall, while the rightmost end is free. External forces acting horizontally on the individual masses are denoted by $u_i$. From Newton's and Hooke's laws we obtain the LTI system
\[
    \dot x = \begin{bmatrix}
        0 & I_N\\ -\calM^{-1} \calK & -\calM^{-1}\calD
    \end{bmatrix} x + \begin{bmatrix}
        0\\ \calM^{-1} E_\mathrm{a}
    \end{bmatrix} u,\qquad y = \begin{bmatrix}
        E_\mathrm{s}^\top & 0
    \end{bmatrix}x
\]
with the mass matrix $\calM = \operatorname{diag} (m_1,\dots, m_N)$, the stiffness matrix $\calK=\mathsf{T}(k_1,\dots,k_N)$
and the damping matrix $\calD=\mathsf{T}(d_1,\dots,d_N)$,
where
\begin{equation*}
    \mathsf{T}(a_1,\dots,a_N) := \begin{bmatrix}
        a_1+a_2 & -a_2\\
        -a_2 & a_2+a_3 & -a_3\\
        & \ddots & \ddots & \ddots\\
        &&-a_{N-1} & a_{N-1} +a_{N} & -a_{N}\\
        &&& -a_{N} & a_{N}\\
    \end{bmatrix}.
\end{equation*}
The matrices $E_\mathrm{a} \in\{0,1\}^{N\times m}$ and $E_\mathrm{s}\in\{0,1\}^{N\times p}$ encode the actuator (external forces) and the sensor placement, respectively, each consisting of columns of the identity matrix $I_N$.

\begin{figure}
    \centering
    \includegraphics[width=0.75\linewidth]{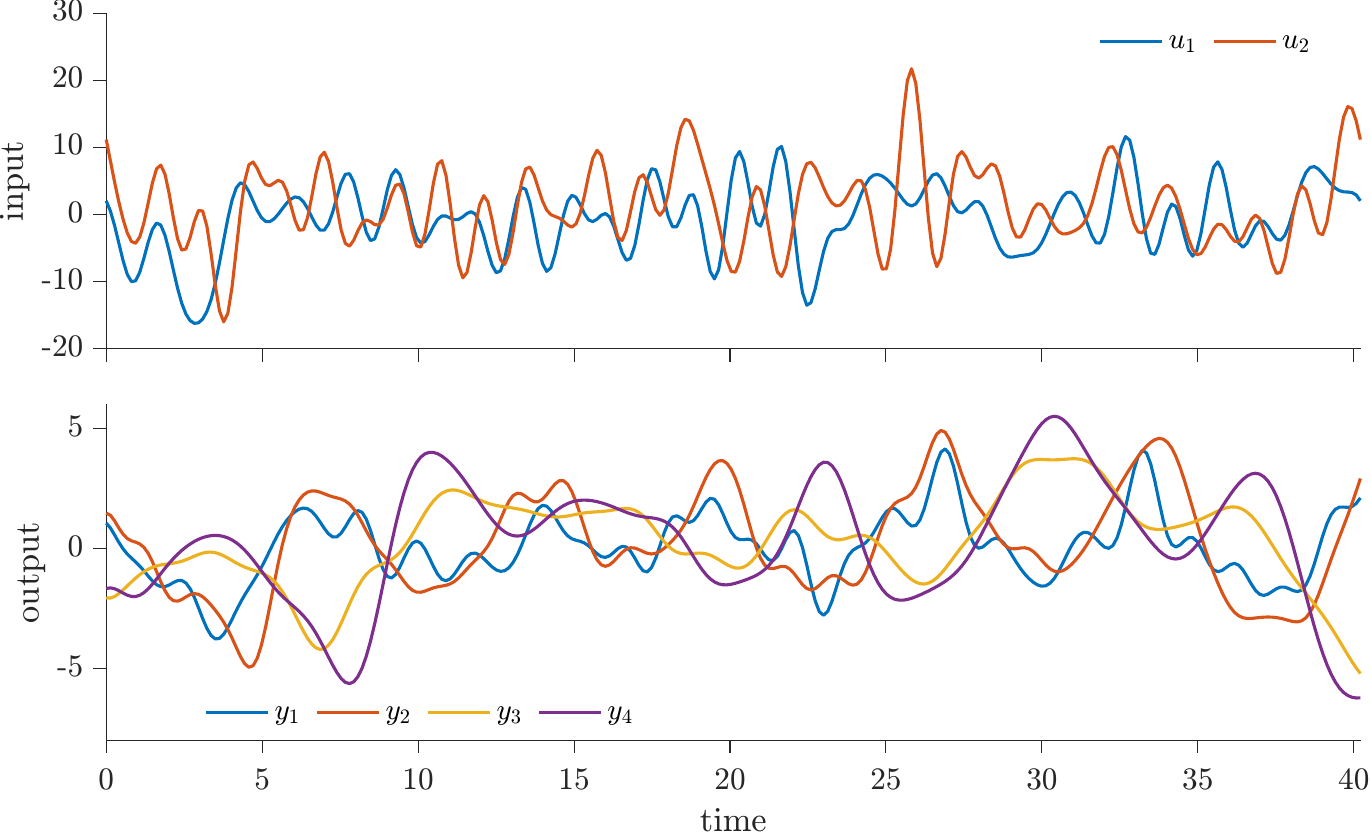}
    \caption{Data trajectory with persistently exciting input signal.}
    \label{fig:pe}
\end{figure}

We choose the chain parameters according to \cref{tab:chain_param}, giving state dimension $n=20$, input dimension $m=2$, output dimension $p=4$ and lag $\lag(\B)=6$. The controllability indices are $\nu_1 = 18$ and $\nu_2=2$. As a persistently exciting input signal of order $L+n$ with $L=\lag(\B)+1=7$, we use a trigonometric polynomial with random coefficients. The output signal is obtained by numerically integrating the state, initialized at a random value. The data trajectory is shown in \cref{fig:pe}.

The resulting Gramian $\Gamma$ has an effective condition
number of $\kappa_+(\Gamma)\approx 1.2128\cdot 10^{13}$.
Preconditioning via Algorithm~3 yields
$\widetilde\Gamma=P\Gamma P^\top$ with
$\kappa_+(\widetilde\Gamma)\approx 2.0991\cdot 10^7$, a reduction by almost
six orders of magnitude, while preserving the numerical rank
$\rank\widetilde\Gamma=Lm+n=34$.

To validate the image representation $M$ obtained by Algorithm~2 together
with \cref{lem:imrep-recovery}, we independently compute a minimal kernel
representation $R$ from the system matrices via the Kalman observability
matrix. We obtain the relative residual error
$$\frac{\lVert R M\rVert}{\lVert R\rVert \lVert M\rVert} \approx 1.153 \times 10^{-12},$$
where $\lVert P\rVert$ for a polynomial matrix $P(s)=\sum_{k=0}^{\deg P} P_k s^k$ is defined as the maximum of its coefficients in the Frobenius norm, i.e.\ $\lVert P\rVert = \max_{0\leq k\leq\deg P}
\lVert P_k\rVert_{\mathrm F}$.
This small residual is consistent with the identity $R(s)M(s)=0$ predicted
by the theory and provides an independent numerical validation of the
data-driven construction.

\begin{figure}
    \centering
    \includegraphics[width=0.75\linewidth]{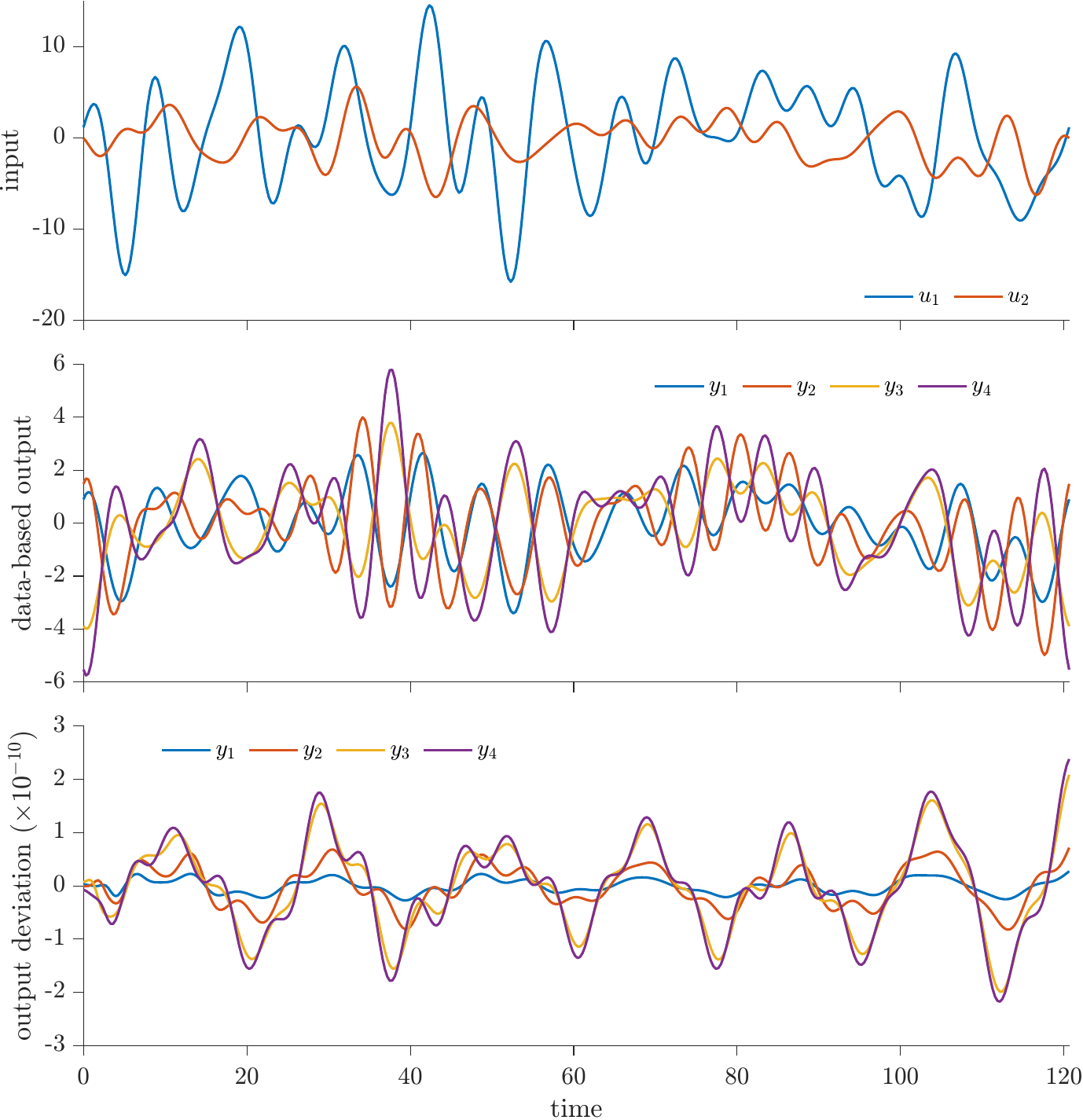}
    \caption{Comparison of data-based and model-based trajectories.}
    \label{fig:iotraj}
\end{figure}

Finally, we generate a latent variable $h$ and calculate the corresponding input--output trajectory via $M$. Since the image representation is not tied to the time interval of the data, we choose a horizon three times as long as that of the persistently exciting input trajectory. For validation, we recover the corresponding initial state and simulate the trajectory directly from the system matrices. \Cref{fig:iotraj} compares the data-based and model-based trajectories, showing the pointwise deviation in each output component.

\begin{table}
\caption{Chain parameters.}
\label{tab:chain_param}
\centering
\begin{tabular}{cccccc}
    %\toprule
    length $N$ & mass $m_i$ & stiffness $k_i$ & damping $d_i$ & actuators at & sensors at\\
    \midrule
     $10$ & $1$ & $5$ & $0.5$ & $1,2$ & $1,3,8,10$\\
    %\bottomrule
\end{tabular}
\end{table}

\bibliographystyle{siamplain}
\bibliography{refs}  
\end{document}